\documentclass[letterpaper, 10 pt, conference]{ieeeconf} 
\usepackage{times} 
\usepackage{amsmath,lipsum} 
\usepackage{amssymb,amsfonts}  
\usepackage{latexsym,theorem,epsfig}
\usepackage{graphicx} 
\usepackage{dsfont}
\usepackage{mathtools}
\usepackage{subfigure}
\usepackage{lscape}
\usepackage{color}
\usepackage{float}
\usepackage{enumerate}
\usepackage[usenames,dvipsnames]{pstricks}
\usepackage{epsfig}
\usepackage{pst-grad} 
\usepackage{pst-plot} 
\usepackage{mathrsfs}
\usepackage{extarrows}
\usepackage[percent]{overpic}

\newtheorem{theorem}{Theorem}

\newtheorem{assumption}{Assumption}

{\theorembodyfont{\rmfamily} }

\allowdisplaybreaks
\IEEEoverridecommandlockouts
\DeclareMathOperator*{\argmin}{arg\,min}
\DeclareMathOperator*{\argmax}{arg\,max}
\title{A Tutorial on Distributed (Non-Bayesian) Learning: \\ Problem, Algorithms and Results}

\author{Angelia Nedi\'{c}, Alex Olshevsky and C\'{e}sar A.\ Uribe
\thanks{Nedi\'{c} is with the ECEE Department, Arizona State University, e-mail: \textit{Angelia.Nedich@asu.edu}. 
Olshevsky is with the ECE Department, Boston University, e-mail: \textit{alexols@bu.edu}. 
Uribe is with the Coordinated Science Laboratory, University of Illinois at Urbana-Champaign,
e-mail: \textit{cauribe2@illinois.edu}.
	This research is supported partially by the National Science Foundation under
	grants CNS 15-44953 and AFOSR FA-95501510394.}  
}
\begin{document}
\maketitle
\begin{abstract}
We overview some results on distributed learning with focus on a family of recently proposed algorithms known as non-Bayesian social learning. 
We consider different approaches to the distributed learning problem and its algorithmic solutions for the case of finitely many hypotheses. 
The original centralized problem is discussed at first, and then followed by a generalization to the distributed setting. 
The results on convergence and convergence rate are presented for both asymptotic and finite time regimes. 
Various extensions are discussed such as those dealing with directed time-varying networks, 
Nesterov's acceleration technique and a continuum sets of hypothesis. 
\end{abstract}

\section{Introduction}

Achieving global behaviors by repeatedly aggregating local information without complete knowledge of the network has been a recent topic of interest~\cite{jad12,rah10,olf06,ala04,zho15}. For example, 
distributed hypothesis testing method that uses belief propagation has been studied in~\cite{ala04}.  
Various extensions to finite capacity channels, packet losses, delayed communications and tracking 
where developed in~\cite{sal06,rah07}. In~\cite{rah10}, the authors proved convergence in probability, asymptotic normality and provided conditions under which the distributed estimation is as good as a centralized one. Later in~\cite{jad12,jad13}, almost sure convergence of non-Bayesian rules based on consensus were shown for static graphs. Other methods to aggregate Bayes estimates in a network have been explored as well \cite{ban14}. The work in~\cite{liu14} extends the results of~\cite{jad12} to time-varying undirected graphs. In~\cite{sha13}, local exponential rates of convergence for undirected gossip-like graphs are studied. The authors in~\cite{lal14b,sha14} proposed a non-Bayesian learning algorithm where a local Bayes' update is followed by a consensus step. In~\cite{lal14b}, convergence result for fixed graphs is provided and large deviation convergence rates are given, proving the existence of a random time after which the beliefs will concentrate exponentially fast. In~\cite{sha14}, similar probabilistic bounds for the rate of convergence are derived and comparisons with the centralized version of the learning rule are provided.

Following the seminal work of Jadbabaie et al.\ in \cite{jad12,mol15,mol12}, there have been many studies of 
non-Bayesian rules for distributed learning. 
Non-Bayesian algorithms involve an aggregation step, usually consisting of a {\it belief aggregation} 
and a {\it Bayesian update} that is based on the locally available data. 
The belief aggregation is typically consisting of a weighted geometric or arithmetic average of beliefs, in which case
the results from consensus literature~\cite{ace08,tsi84,jad03,ned13,ols14} are exploited, while Bayesian update step is based on
the standard Bayesian learning approach~\cite{ace11,mos14}.
 
Several variants of non-Bayesian approach have been proposed and have been shown 
to produce consistent estimates, with provable asymptotic and non-asymptotic convergence rates 
for a general class of distributed algorithms. The main body of work is focused on the case of finitely many hypotheses.
The established results include
asymptotic convergence rate analysis~\cite{sha13,lal14,qip11,qip15,sha15,rah15,sah16} and 
non-asymptotic convergence rate bounds~\cite{sha14,ned15,lal14b}, 
time-varying directed graphs \cite{ned15}, continuum set of hypotheses \cite{nedic2016b}, weakly connected graphs \cite{sal16}, bisection search algorithm \cite{tsi16}, and transmission node failures ~\cite{su16,su16b,su16c}. 

In this paper, we overview a subset of recent studies on distributed (non-Bayesian) learning algorithms. 
To present a concise introduction to the topic, we start by presenting ideas from centralized learning and, then, 
transition to the most recent developments in the distributed setting. 
This tutorial is by no means extensive and the interested reader may like to look into 
the references for a more complete exposition of certain aspects.

This tutorial is organized as follows. Section \ref{problem} presents a general introduction to the distributed learning problem. We highlight the main assumptions and how they can be weakened for more general results. 
Section \ref{algos} provides an overview of the centralized non-Bayesian learning problem and 
describes some initial generalizations to the distributed setting (known as social learning). 
Moreover, convergence results as well as (non-)asymptotic convergence rate estimates are provided. 
Section \ref{extension} discusses some generalizations aimed at improving the convergence rate estimates
(in terms of their dependency on the number of agents), 
dealing with time-varying directed graphs, and learning with a continuum sets of hypotheses. 
Finally, some conclusions are presented in Section \ref{conclusions}.
 
 \textbf{\textit{Notation:}} 
 The inner product of two vectors $x$ and $y$ is denoted by $\langle x,y\rangle$.
 We write $\left[A\right]_{ij}$ or $A_{ij}$ to denote the entry of a matrix $A$ in 
 the $i$-th row and $j$-th column. We write $A'$ for the transpose of a matrix $A$ and 
 $x'$ for the transpose of a vector $x$. 
 A matrix is said to be stochastic if its entries are nonnegative, and the sum of the entries in every row is equal to 1.
 A stochastic matrix $A$ whose transpose $A'$ is also stochastic is said to be doubly stochastic.
 We use $I$ for the identity matrix, where its size will be inferred from the context. 
 We write $e_i$ to denote the vector with all zero entries except for its $i$-th entry which is equal to $1$. 
 In general, when referring agents we will use superscripts with the letter $i$ or $j$, while when referring to a time instant 
 we will use subscripts and the letter $k$. 
 
 We write $|\Theta|$ to denote the cardinality of a set $\Theta$,
and $\Delta_\Theta$ for a probability measure over the set $\Theta$. 
Upper case letters represent random variables (e.g. $X_k$) 
with their corresponding lower case letters as their realizations (e.g. $x_k$). Notation $\mathbb{E}_X$ is reserved for expectation with respect to a random variable $X$. 
We denote the Kullback-Liebler (KL) divergence between two probability distributions $p$ and $q$ 
 with a common support set by 
 $D_{KL}\left(p \| q\right)$. In particular, when the distributions $p$ and $q$ have a countable (or finite) support set,
 their KL-divergence is given by
 \begin{align*}
 	D_{KL}\left(p \| q\right)  = \sum_{i=1}^\infty p_i \log\left(\frac{p_i}{q_i}\right).
	 \end{align*} 
The definition of the KL-divergence for general measures $p$ and $q$ on a given set is a bit more involved;
it can be found, for example, in~\cite{bernardo2001bayesian}, page 111.
	 
\section{Problem Statement}\label{problem}
\begin{figure}[t]
	\centering
	\begin{overpic}[width=0.25\textwidth]{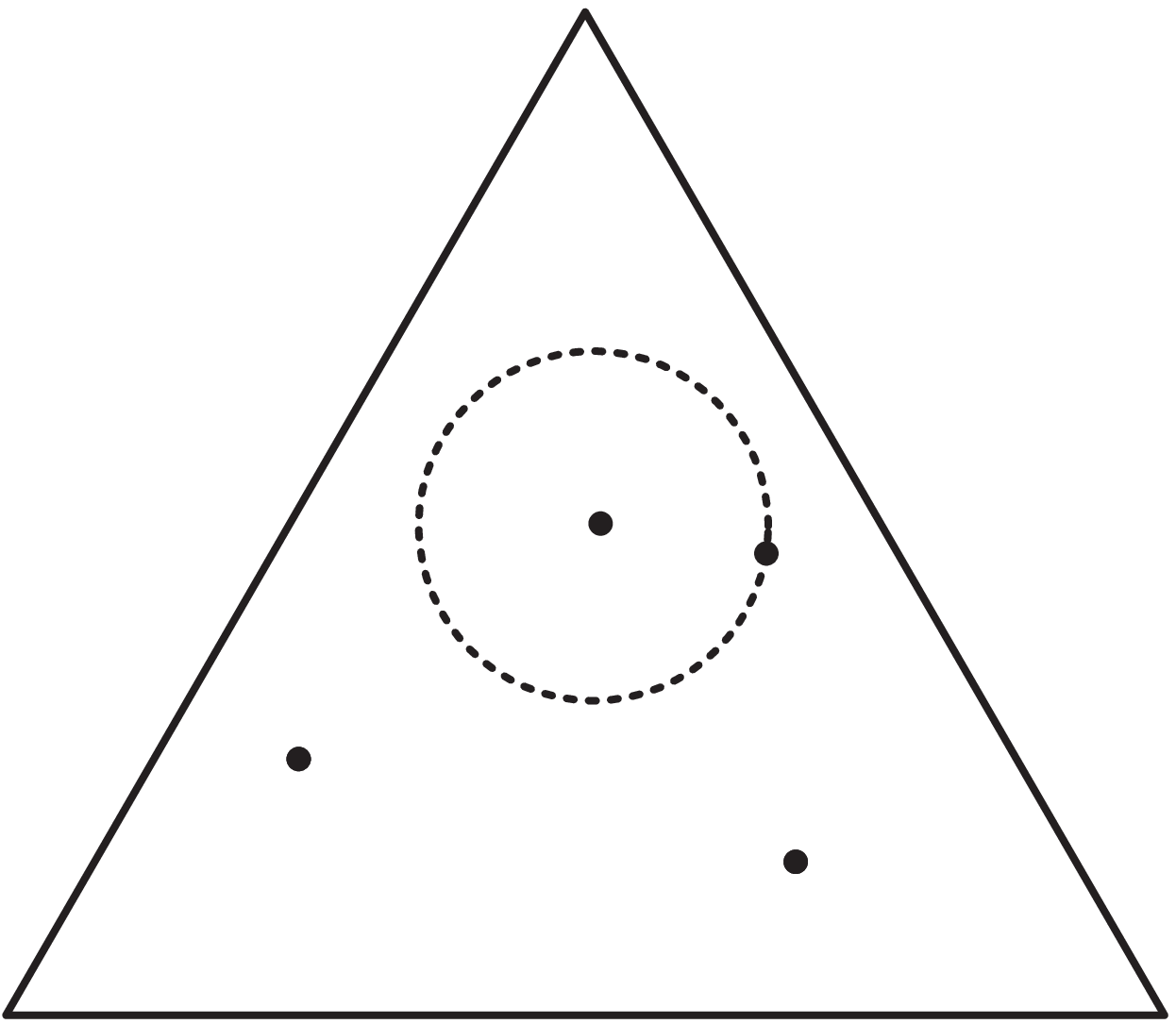}
		\put(50,46){{\small $\boldsymbol{f}$}}
		\put(65,33){{\small $\boldsymbol{\ell}(\cdot|\theta^*)$}}
		\put(27,23){{\small $\boldsymbol{\ell}(\cdot|\theta_1)$}}
		\put(65,16){{\small $\boldsymbol{\ell}(\cdot|\theta_2)$}}
	\end{overpic} 
	\caption{Geometric interpretation of the learning objective. The triangle represents the simplex composed of all agents' probability distributions. The observations of the agents are generated according to a joint probability distribution 
		$\boldsymbol{f}\left(\cdot\right)$. The joint distribution for the agent observations is parametrized by $\theta$. 
		The agent goal is to learn a hypothesis that best describes their observations, which corresponds to the distribution 
		$\ell(\cdot|\theta^*)$ (the closest to the distribution $\boldsymbol{f}\left(\cdot\right)$).}
	\label{triangle}
\end{figure}

Consider a group of $n$ agents, indexed by $1,2,\ldots,n$, 
each having conditionally independent observations of a random process at discrete time steps $k=1, 2, 3, \ldots$. 
Specifically,  agent $i$ observes the random variables $S_1^i, S_2^i, \ldots,$ which are i.i.d.\ in time 
and distributed according  to an unknown probability distribution $f^i$. 
The set of possible outcomes of the random variables $S_k^i$ is a finite set which we will denote by
${\mathcal{S}}^i$. For convenience, we stack up all the $S_k^i$ into a vector denoted as $\mathbf{S}_k$. 
Then, $\mathbf{S}_k$ is an i.i.d.\ vector taking values in 
${\boldsymbol{\mathcal{S}} = \prod_{i=1}^n \mathcal{S}^i}$ and distributed as ${\boldsymbol{f}=\prod_{i=1}^n f^i}$. Furthermore, each agent $i$ has a family of probability distributions $\{\ell^i(\cdot | \theta)\}$ 
parametrized by $\theta \in\Theta$, where $\Theta$ is a set of parameters. 
One can think of $\Theta$ as a set of hypotheses and $\ell^i( \cdot | \theta)$ as the probability distribution that would be seen by agent $i$ if hypothesis $\theta$ were true. In general it is not required that there exists
$\theta \in \Theta$ with $\ell^i(\cdot | \theta) = f^i$ for all $i=1, \ldots, n$; in other words, there may not be a hypothesis which matches the observations made by the nodes. Rather, the objective of all agents is to 
agree on a subset of $\Theta$ that fits all the observations in the network best. 
Formally, this setup describes the scenario where the group of agents collectively tries to solve the following optimization problem
\begin{align}\label{opt_problem}
\min_{\theta \in \Theta} F(\theta) & \triangleq D_{KL}\left(\boldsymbol{f}\|\boldsymbol{\ell}\left(\cdot | \theta\right)\right) \ \ \   \\
& = \sum\limits_{i=1}^n D_{KL}\left(f^i\|\ell^i\left(\cdot|\theta\right)\right), \nonumber
\end{align}
where $D_{KL}\left(f^i\|\ell^i\left(\cdot|\theta\right)\right)$ is the Kullback-Leibler divergence between the distribution of $S_k^i$ and the distribution $\ell^i( \cdot | \theta)$ that would have been seen by agent $i$ if hypothesis $\theta$ were correct. 
The distributions $f^i$'s are not known, therefore, the agents want to ``learn" the solution to this optimization problem based on their local observations and some local interactions. See Figure~\ref{triangle} for an illustration of the problem. 

The agents interact over a sequence of directed communication graphs 
$\mathcal{G}_k = \{V,E_k\}$, where $V = \{1,2,\ldots,n\}$ is the set of agents (where each agent is viewed as a node), and
$E_k$ is the set of edges where $(j,i) \in E_k$ if agent $j$ can communicate with agent $i$ at time instant $k$. 
Specifically, the agents communicate with each other by sharing 
their beliefs about the hypotheses set, denoted as $\mu^i_k$, which is a probability distribution over the hypothesis set $\Theta$. 
In the forthcoming discussion, we will consider the cases where the graphs $\mathcal{G}_k$ can be static and may be undirected.
We will clearly specify the assumptions made on the graphs.

The hypothesis set $\Theta$ can be a finite, countable or continuum set, 
which will be self-evident from expressions used in Bayes' update relation. 

\section{Algorithms}\label{algos}
In this section, we describe  some of the algorithms 
that have been proposed for the distributed non-Bayesian learning problem. 
Different algorithms and results exist due to the use of different communication networks and protocols for information exchange. Moreover, the variety in the algorithms is also due to the order in which the local information updates and neighbor beliefs aggregation updates are performed. 

We will start by considering Bayes' update for a case of a single agent, i.e., centralized case. Furthermore,
initially, for simplicity of exposition we will assume there exists a single $\theta^*$ 
that minimizes problem \eqref{opt_problem} corresponding to a single agent case. 
In this case, updating the beliefs to account for a set of observations that lead to a posterior belief follows the Bayes'rule. Specifically, having a belief $\mu_k$ and a new observation $s_{k+1}$ at time $k$, the agent updates its belief as follows 
(see~\cite{eps10}): 
\begin{align}\label{bayes}
\mu_{k+1} & = BU(\mu_k; s_{k+1}), 
\end{align}
where $BU(\mu_k; s_{k+1}) $ denotes the Bayesian update of the belief $\mu_k$ given a new observation $s_{k+1}$, i.e.,
\begin{align*}
	\mu_{k+1}(\theta) & = \frac{\mu_k(\theta)\ell(s_{k+1}|\theta
		)}{\int_{\Theta}d\mu_k(\theta)\ell(s_{k+1}|\theta)} \\
	& \propto \mu_k(\theta)\ell(s_{k+1}|\theta)\qquad\hbox{for all }\theta\in\Theta,
\end{align*}
where the symbol $\propto$ stands for positively proportional quantities 
(the proportionality constant here is the normalization factor needed to obtain a probability distribution).

In \cite{eps10,eps2008}, the concepts of over-reaction and under-reaction to local observations are introduced,
where the belief update rule is given by
\begin{align}\label{under_over}
\mu_{k+1} & = (1- \gamma_{k+1})BU(\mu_k; s_{k+1}) +\gamma_{k+1}\mu_k,
\end{align}
where $\gamma_{k+1} \leq 1$. 
When $\gamma_{k+1} = 0$, 
algorithm \eqref{under_over} reduces to Bayesian learning in \eqref{bayes}. 
When $\gamma_{k+1} > 0$, a relative importance is given to the prior, 
whereas, for $\gamma_{k+1}<0$ the updates over-react to observations. 
The authors in \cite{eps10} showed that update rules of the form \eqref{under_over} converge to the correct parameter value 
$\theta^*$ in the almost sure sense whenever $0 \leq \gamma_{k+1} \leq 1$ and $\gamma_{k+1}$ is measurable. 
If $\gamma_{k+1} < 0$ or if $\gamma_{k+1}$ is not measurable, then there is an incorrect parameter $\bar{\theta}$ 
to which convergence can happen with a positive probability. Thus, as long as there is a constant flow of new information and 
the agent takes its personal signals into account in a Bayesian manner, 
the resulting learning process asymptotically coincides with Bayesian learning.

The seminal work of \cite{jad12} has introduced a social learning approach to non-Bayesian learning, where different agents receive different observations and use a DeGroot-style update to incorporate the views of their neighbors
\begin{align}\label{nonB_learning}
\mu_{k+1}^i & = a_{ii}BU(\mu_k; s_{k+1}) +\sum_{j=1}^{n}a_{ij}\mu_k^j,
\end{align}
where $a_{ij}\ge0$ are the weights taking positive values on the links $(j,i)$ in a static graph 
(i.e., $\mathcal{G}_k=\mathcal{G}$ for all $k$) and satisfying $\sum_{j=1}^n a_{ij}=1$ for all $i$.  
In~\cite{jad12}, it has been shown that, when the underlying social network $\mathcal{G}$ 
is strongly connected, every $a_{ii}>0$, and at least one agent has a positive prior belief on the true parameter
(i.e., $\mu^i(\theta^*) >0$ for some $i$), then the beliefs 
generated by algorithm~\eqref{nonB_learning} results in all agents forecasts converging to the 
correct one with probability one.

A connection between non-Bayesian learning and optimization theory were pointed out in \cite{sha13}, where a distributed learning algorithm has been proposed that is based on a maximum likelihood analysis of the estimation problem 
and Nesterov's dual averaging algorithm~\cite{nes09}. 
Finding the true state of the world was described as the following optimization problem
\begin{align}\label{opt1}
\max_{\mu \in \Delta_\Theta} 
f(\mu) , \qquad f(\mu)\triangleq \left\langle\mu,\sum\limits_{i=1}^{n} \mathbb{E}_{S^i}\left[\log \ell^i (S^i|\theta)\right]
\right\rangle, 
\end{align}
or equivalently
\begin{align}\label{opt2}
\theta^* & = \argmax_{\theta \in \Theta} \left\lbrace \mathbb{E}_{S} \left[\ell (S|\theta)\right] \right\rbrace.  
\end{align}
Applying a regularized dual averaging algorithm to the optimization problem \eqref{opt2}, 
one obtains a sequence $\{\mu_k,z_k\}_{k=0}^{\infty}$, where 
\begin{subequations}\label{dual_nesterov}
	\begin{align}
		z_{k+1} & = z_k +g_k,\\
		\mu_{k+1} & = \argmin_{x\in \Delta_\Theta} \left\lbrace \left\langle z_{k+1},x\right\rangle  + \frac{1}{\alpha_k} \psi (x)\right\rbrace, 
	\end{align}
\end{subequations}
with $\mu_k \in \Delta_\Theta$, $z_k \in \mathbb{R}^m$, $\{\alpha_k\}_{k=0}^{\infty}$ being a sequence of non-increasing step-sizes, $g_k = \log \ell(s_k|\theta)$ and $\psi (x)$ a \textit{proximal function}.

Specifically for the centralized case with the Kullback-Liebler divergence as proximal function, the algorithm in \eqref{dual_nesterov} has an explicit closed form solution which coincides with \eqref{bayes}.

In the distributed setting in \cite{sha13}, for an undirected and static graph $\mathcal{G}$,
the randomized gossip interactions were considered, 
where an agent $i$ ``wakes-up" according to a Poisson clock and communicates with a randomly selected agent $j$. Both agents average their accumulated observations and add their most recent stochastic gradient, resulting in the update of the form:
\begin{subequations}\label{dist_dual_nesterov}
	\begin{align}
	z^i_{k+1} & = \frac{1}{2}(z_k^i + z_k^j) + g^i_k, \\
	z^j_{k+1} & = \frac{1}{2}(z_k^i + z_k^j) + g^j_k, \\
	\mu_{k+1}^{(i,j)} & = \argmin_{x\in \Delta_\Theta} \left\lbrace \langle z_{k+1}^{(i,j)},x\rangle  
	+ \frac{1}{\alpha_k} \psi (x)\right\rbrace, 
	\end{align}
\end{subequations}
while the other agents in the system do not update.

Letting $z^i_0 = 0$ for all agents $i$ and using the Kullback-Liebler divergence as a proximal function, 
the update rule of the form~\eqref{dist_dual_nesterov} has a closed form solution given by
\begin{align}
	\mu^i_{k}(\theta) & \propto \mu^i_0(\theta) \exp \left( k \phi^i_k \left(\theta\right) \right) \qquad\hbox{for all }i,
\end{align}
where 
\begin{align*}
\phi^i_k = \frac{1}{k} \sum\limits_{\tau=0}^{k-1} \sum\limits_{j=1}^{n} \prod\limits_{\rho=1}^{k-1-\tau}\left[ A_{k-\rho}\right]_{ij} \log \ell^j(s^j_\tau|\theta) 
\end{align*}
and $
A_k  = I - \frac{\left( e_{i_k} - e_{j_k}\right)\left( e_{i_k} - e_{j_k}\right)' }{2},
$
with $i_k$ and $j_k$ being the agents involved in the random gossip communication at time $k$ (or alternatively, the link 
$\{i_k,j_k\}\in E$ being randomly activated in the graph $\mathcal{G}=(V,E)$).

The update rule in~\eqref{dist_dual_nesterov} 
involves a form of geometric average of beliefs instead of the linear aggregations of beliefs as in~\eqref{nonB_learning}. 
Weak convergence is proven under the connectivity assumption of the interacting graph $\mathcal{G}$, i.e.,
\begin{align*}
	\mu^i_k(\theta^*) \to 1  \text{ as } k \to \infty \text{ for all } i.
\end{align*}

Additionally, in \cite{sha13}, the convergence rate results for the estimation process are provided. 
An asymptotic rate is derived that guarantees that for sufficiently large time scales the beliefs will concentrate around the true hypothesis exponentially fast. The rate at which this happens is proportional to the distance (in the sense of the KL divergence) between the true hypothesis and the second best option, i.e., with probability $1 - \delta(\epsilon, k )$ for sufficiently large $k$ it holds that
\begin{align*}
|\mu^i_k(\theta^*) -1 | & \leq \mathcal{K} \exp((-D + \epsilon) k),
\end{align*}
where $\mathcal{K}$ is a constant and $
D = \argmin_{\theta \in \Theta \setminus \{\theta^*\}} \frac{1}{n} F(\theta).
$

Similar asymptotic rates using large scale deviation theory were derived in \cite{lal14b} for a directed static graph 
but for a different algorithm.
Specifically, in \cite{lal14b}, an explicit belief update rule is considered
where local Bayesian updates are aggregated via geometric averages of the form: 
\begin{align}\label{non_bayes_distributed_e}
\mu_{k+1}^i\left(\theta\right) & \propto  \prod_{j=1}^n BU\left(\mu_k^j\left(\theta\right);s_{k+1}^j\right)^{A_{ij}}.
\end{align} 
Under assumptions of strong connectivity, positive prior beliefs and existence of unique correct models, an exponential convergence rate of the beliefs to the correct hypothesis has been shown and an asymptotic convergence rate is provided (see Theorem 1 of \cite{lal14b}).

In recent works \cite{sha14,ned14}, non-asymptotic convergence rates for a variety of distributed non-Bayesian learning algorithms
have been established. 
In \cite{sha14}, the algorithm in \eqref{dist_dual_nesterov} has been considered for the case of 
(non-random) agent interaction over a general static connected graphs. In particular, the following relations have been shown to hold
\begin{align}\label{shahin}
& \frac{1}{\eta} \log \|\mu_k^i - e^*\|_{TV} \cr
& \leq - \argmin_{\theta \in \Theta \setminus\{\theta^*\} }
k \sum_{j=1}^{n} \pi^j D_{KL}(\ell^i(\cdot|\theta^*)\|\ell^i(\cdot|\theta)) \nonumber \\
& + \sqrt{2(\log \alpha )^2 k \log \frac{|\Theta|}{\delta} }+ \frac{8 \log \alpha\log n}{1 -\lambda_{max}(A)} 
+ \frac{\log |\Theta|}{\eta}\qquad
\end{align}
with probability $1-\delta$, where $\delta>0$ is arbitrarily small. Here, $\|a- b\|_{TV}$ denotes the total variation between vectors $a$ and $b$,  $e^*$ is a probability vector with an  entry $1$ in the position corresponding to hypothesis $\theta^*$, $|\Theta|$ denotes the size of the hypotheses set, $\eta>0$ is a step-size and $\alpha>0$ is a lower bound on the probability mass distribution in the likelihood models. 
The vector $\pi$ denotes the stationary distribution of the corresponding Markov chain whose transition probability distribution is the interaction matrix $A$ (in other words the vector $\pi$ is the left eigenvector associated 
with the eigenvalue 1 of the matrix $A$).

The non-asymptotic probabilistic bound in \eqref{shahin} shows the concentration of the beliefs around the true state of the world as an exponentially fast process with a transient time related to the matrix properties and the desired accuracy level. 
The bound holds for a connected graph and stochastic weight matrix $A$, and 
the exponential concentration rate depends explicitly on the left-eigenvector $\pi$ associated with eigenvalue 1 of the matrix $A$.

An independent simultaneous work \cite{ned14,ned15}
also has developed non-asymptotic bounds for distributed non-Bayesian learning for time-varying graphs
and for different algorithms. 
The belief update rules in \cite{ned14,ned15} are based on mirror descent algorithm as applied to the learning problem in a distributed setting. 
The resulting update rule has the following form:
\begin{align}\label{non_bayes_distributed2}
\mu_{k+1}^i\left(\theta\right) & \propto   BU\left(\prod_{j=1}^n\mu_k^j\left(\theta\right)^{\left[A_k\right]_{ij}};s_{k+1}^i\right).
\end{align} 
Algorithm \eqref{non_bayes_distributed2} is applicable to time-varying graphs, as indicated by the use of time varying weight matrices $A_k$ that are compliant with the graphs' structure. 
In particular, the following assumption is imposed on the graph sequence $\left\{\mathcal{G}_k\right\}$ and the matrix sequence
$\{A_k\}$.

\begin{assumption}\label{assum:graph}
Assume that each graph $\mathcal{G}_k$ is undirected and has no self-loops (i.e., $\{i,i\}\not \in E_k$ for all $i$ and all~$k$). 
Moreover,
	let the graph sequence $\{ \mathcal{G}_k \}$ and the matrix sequence $\{A_k\}$ satisfy the following conditions:
	\begin{enumerate}[(a)]
		\item $A_k$ is doubly-stochastic for every $k$, with $\left[A_k\right]_{ij} > 0$ if $\{i,j\}\in E_k$ and 
		 $[A_k]_{ij}=0$ for $\{i,j\} \notin E_k$.
		\item Each $A_k$ has positive diagonal entries, i.e., $\left[A_k\right]_{ii}>0$ for all $i=1,\ldots,n$ and all $k\ge0$.
		\item There exist a uniform lower bound $\eta>0$ on positive entries in $A_k$, i.e.,
		$\left[A_k\right]_{ij} \ge \eta$ if $\left[A_k\right]_{ij}>0$.
		\item The graph sequence $\left\{\mathcal{G}_k\right\}$ is $B$-connected, 
		i.e.,
		there is an integer $B\ge 1$ such that the graph 
		$\left\{V,\bigcup_{i=kB}^{\left(k+1\right)B-1}E_i\right\}$ is connected for all $k \geq 0$.
	\end{enumerate}
\end{assumption}

We are now considering the learning problem in~\eqref{opt_problem}, where the hypothesis set $\Theta$ is finite.
We let $\Theta^*$ denote the set of optimal solutions, and note that this set is nonempty.
In this setting, the following assumption ensures that the learning process will identify correct hypothesis.
In particular, the assumption is for the general case when a unique true state $\theta^*$ of the underlying process does not exist
(implying that $\Theta^*$ is not a singleton).

\begin{assumption}\label{assum:init}
	For all agents $i=1,\ldots,n$,
	\begin{enumerate}[(a)]
		\item
		There is a nonempty set $\Theta^{*i}\subseteq \Theta^*$ such that $\mu_{0}^i\left(\theta\right)>0$ 
		for all  $\theta \in \Theta^{*i}$. Furthermore, the intersection set $\cap_{i=1}^n\Theta^{*i}$ is nonempty.
		\item 
		There exists an $\alpha >0$ such that if $f^i\left(s^i\right)>0$ then  
		$\ell^i\left(s^i | \theta \right)  > \alpha$ for all $\theta \in \Theta$.
	\end{enumerate}
\end{assumption}
With the two assumptions above we can state the main result in \cite{ned15}.
\vspace{-1em}
\begin{theorem}\label{teo2} 
	Let Assumptions \ref{assum:graph} and \ref{assum:init} hold, and let $\rho\in(0,1)$. The update rule of Eq.~(\ref{non_bayes_distributed2}) has the following property:
	there is an integer $\boldsymbol{N}(\rho)$ such that, with probability $1 -\rho$, 
	for all $k\ge \boldsymbol{N}(\rho)$ and for all $\theta\notin\Theta^*$, we have
	$$\mu_{k}^i\left(\theta \right) \leq \exp\left( -\frac{k}{2}\gamma_2+ \gamma_1^i\right)
	\quad\hbox{for all } i = 1, \ldots, n,$$ 
	where $
		\boldsymbol{N}(\rho)
		\triangleq \left\lceil\frac{8 \left( \log\left(\alpha \right)\right) ^2  \log\left(\frac{1}{\rho} \right) }{\gamma_2^2} + 1 \right\rceil,
	$
	{\small
	\begin{align*}
		\gamma_1^i & \triangleq \max_{\substack{\theta\notin \Theta^*\\ \theta^*\in \hat{\Theta}^*} } \left\{
		\log \frac{\mu_0^i(\theta)}{\mu_0^i(\theta^*)}
		+ {\frac{8 \log n}{1 - \lambda} \log\frac{1}{\alpha}} + F^i\left(\theta\right)-F^i\left(\theta^*\right) \right\},\cr
		\gamma_2 &\triangleq \frac{1}{n}\, \min_ {\theta\notin\Theta^*}
		\left(F(\theta)-F(\theta^*)\right),
	\end{align*}
}
	\hskip-0.5pc where $F^i(\theta)$ is a local learning objective of agent $i$ given by
	\[
	F^i(\theta) = D_{KL}\left(f^i\|\ell^i\left(\cdot|\theta\right)\right),\]
	while $\alpha$ is a constant from Assumption~\ref{assum:init}(b), $\eta$ from Assumption~\ref{assum:graph}(c)
	and $\lambda$ is given by $
		\lambda = \left(1-\frac{\eta}{4n^2} \right)^{\frac{1}{B}}.
	$
	If each  $A_k$ is the 
	lazy Metropolis matrix associated with $\mathcal{G}_k$ and $B=1$, then
	$
		\lambda = 1- \frac{1}{\mathcal{O}(n^2) }.
	$
\end{theorem} 

Theorem~\ref{teo2} states that, with a high probability and after some time that is sufficiently long (as captured by $N(\rho)$), 
the belief of each agent on any hypothesis outside the optimal set decays at a network-independent rate. 
This rate scales with the constant $\gamma_2$, which is the
average Kullback-Leibler divergence to the next best hypothesis. However, there is a transient due to the $\gamma_1^i$ term (since the bound of 
Theorem \ref{teo2} is not even below $1$ until $k \geq 2 \gamma_1^i/\gamma_2$), and the size of this transient depends on the network and the
number of nodes through the constant $\lambda$.

We note that the transient time for each agent $i$ 
is affected by the discrepancy in the initial beliefs on the correct hypotheses 
(those in the set $\Theta^*$), as captured by the term
\[\log \frac{\mu_0^i(\theta)}{\mu_0^i(\theta^*)}\]
in the expression for $\gamma_1^i$ in Theorem~\ref{teo2}.
We note that, if agent $i$ uses a uniform initial belief, i.e., $\mu_0^i(\theta)=1/|\Theta|$ for all $\theta \in\Theta$,
then this term would be 0 for all $\theta$ and, consequently, it will not contribute to the transient time $\gamma_1^i$.
Thus, the transient time has a dependence on the initial beliefs that is intuitively plausible.
Moreover, if agent $i$ were to start with a good initial belief $\mu_0^i$, i.e., a belief such that 
\[\log \frac{\mu_0^i(\theta)}{\mu_0^i(\theta^*)}<0 \qquad\hbox{for all }\theta\not\in\Theta^*,\]
then the corresponding transient time $\gamma_1^i$ would be smaller, which is also to be expected.

\subsection{Connection with Distributed Stochastic Mirror Descent}

To make this connection simple, we will keep the assumption that the hypothesis set $\Theta$ is a finite set.
Then, we  can observe that 
the optimization problem in Eq.~\eqref{opt_problem} is equivalent to the following problem:
\begin{align*}
\min_{\pi \in \Delta_\Theta } \mathbb{E}_{\pi}  \sum\limits_{i=1}^n D_{KL}\left(f^i\|\ell^i\right)
& = \min_{\pi \in \Delta_\Theta } \sum\limits_{i=1}^n \mathbb{E}_{\pi} \mathbb{E}_{f^i} [- \log \ell^i].
\end{align*} 
The expectations in the preceding relation can exchange the order, so the problem in Eq.~\eqref{opt_problem} 
is equivalent to the following one:
\begin{align}\label{main2}
\min_{\pi \in \Delta_{\Theta} } \sum\limits_{i=1}^n \mathbb{E}_{f^i} \mathbb{E}_{\pi} [- \log \ell^i].
\end{align}
The difficulty in evaluating the objective function in Eq.~\eqref{main2} (even in the case of a single agent)
lies in the fact that the distributions $f^i$ are unknown. 
A generic approach to solving such problems is the class of stochastic approximation methods, where the objective is minimized by constructing a sequence of gradient-based iterates where the true gradient of the objective 
(which is not available) is replaced with a gradient sample that is available at the given update time.
A particular method that is relevant here is the stochastic mirror-descent method which would solve the problem 
in Eq.~\eqref{main2}, in a centralized fashion, by constructing a sequence $\{x_k\}$, as follows:
\begin{align}\label{central}
x_{k+1} & =  \argmin_{y \in X} \left\lbrace \langle g_{k}, y\rangle 
+ \frac{1}{\alpha_{k}} D_w(y,x_{k})\right\rbrace 
\end{align}
where $g_k$ is a noisy realization of the gradient of the objective function 
in Eq.~\eqref{main2} and $D_w(y,x)$ is a Bregman distance function associated with a distance-generating 
function $w$, and $\alpha_k>0$ is the step-size. 
If we take $w(t) = t \log t$ as the distance-generating function, then the corresponding Bregman distance is 
the Kullback-Leiblier divergence $D_{KL}$. 
Let us note that 
this specific selection of Bregman divergence was previously studied in \cite{bec03}, where the entropic mirror descent algorithm was proposed. Thus, in this case, 
{\it the update rule in Eq.~\eqref{protocol_cont} corresponds to
	a distributed implementation of the stochastic mirror descent algorithm} in~\eqref{central},
	where $D_w(y,x)=D_{KL}(y,x)$, $X=\Delta_\Theta$, and the stepsize is fixed, i.e.,  $\alpha_k =1$ for all $k$.
	
The update rule in Eq. \eqref{non_bayes_distributed2} {defines a probability measure $\mu_k^i$over the set $\Theta$} which coincides with the iterate update of the distributed stochastic mirror descent algorithm applied to the optimization problem in 
Eq.~\eqref{opt_problem}, i.e.,	
	
	{\small	
	\begin{align}\label{stoc_mirror}
	\mu_{k}^i = \argmin_{\pi \in \Delta_{\Theta} } \left\lbrace 
	\mathbb{E}_{\pi} [\textendash\log \ell^i(s_{k}^i|\cdot)] + \sum\limits_{j=1}^{n} \left[ A_k\right] _{ij} D_{KL}(\pi\|\mu_{k-1}^j)\right\rbrace. 
	\end{align}
	}

\section{Extensions}\label{extension}

\subsection{Fast Rates with Nesterov's Acceleration}

For static undirected graphs, the authors in \cite{ned15} proposed an update rule with one-step memory as follows:
\begin{align}\label{linear_bayes}
\mu_{k+1}^i\left(\theta\right) & 
\propto 
\frac{ \prod\limits_{j=1}^{n} \mu_k^j\left(\theta\right)^{\left(1+\sigma\right)
		\bar{A}_{ij} } \ell^i \left(s^i_{k+1}|\theta\right) }{\prod\limits_{j=1}^{n} \left( \mu_{k-1}^j\left(\theta\right)\ell^j \left(s^{j}_{k}|\theta\right)\right) ^{\sigma\bar{A}_{ij} }},
\end{align}
where $\sigma$ a constant to be set later. 
This update rule is based on an accelerated algorithm for computing network aggregate values as given in~\cite{ols14}, which 
has the convergence rate a factor of $n$ faster than the previous rate results (in terms of the factor that governs the exponential decay). 

For the algorithm in~\eqref{linear_bayes} we impose the following assumption.
\begin{assumption}\label{assum_linear}
	The graph sequence $\{ {\cal G}_k\}$ is static {(i.e. $\mathcal{G}_k = \mathcal{G}$ for all $k$)} and undirected and the weight matrix
	$\bar{A}$ is a lazy Metropolis matrix, defined by
	\begin{align*}
	\bar{A} = \frac{1}{2}I + \frac{1}{2}A,
	\end{align*}
	where $A$ is the \textit{Metropolis matrix}, which is the unique stochastic matrix whose off-diagonal entries satisfy
	\begin{align*}
	A_{ij} = \left\{  \begin{array}{l l}
	\frac{1}{\max \left\{ d^i+1,d^j+1 \right\} } & \quad \text{if $\{i,j\} \in E$, }\\
	0 & \quad \text{if $\{i,j\} \notin E$,}
	\end{array} \right. 
	\end{align*}
	with $d^i$ being the degree of the node $i$ (i.e., the number of neighbors of $i$ in the graph).
\end{assumption}

The next theorem provides a convergence rate bound for the beliefs generated by algorithm~\eqref{linear_bayes}.
In particular, it shows the rate at which the beliefs dissipate the mass placed on wrong (non-optimal) hypotheses.

\begin{theorem}\label{thm_lineal} 
	Let Assumptions~\ref{assum_linear} and~\ref{assum:init} hold and let $\rho\in(0,1)$. Furthermore let $U \geq n$ and let $\sigma = 1 - 2/(9U+1)$.	Then, the update rule of Eq.~\eqref{linear_bayes} with this $\sigma$, {uniform} initial condition $\mu^i_{-1}\left(\theta\right)=\mu^i_{0}\left(\theta\right)$ and $\beta_{-1}^i$ fixed to zero, has the following property:
	there is an integer $\boldsymbol{N}(\rho)$ such that,
	with probability $1 -\rho$, 
	for all $k\ge \boldsymbol{N}(\rho)$ and for all $\theta\notin\Theta^*$,
	there holds
	$$\mu_{k}^i\left(\theta \right) \leq \exp\left( -\frac{k}{2}\gamma_2+ \gamma_1\right)
	\quad\hbox{for all  }i = 1, \ldots, n,$$ 
	where $
	\boldsymbol{N}(\rho)
	\triangleq \left\lceil\frac{72 \left( \log\left(\alpha \right)\right) ^2 n 
	\log\left(\frac{1}{\rho} \right) }{\gamma_2^2} \right\rceil,
$
	\begin{align*}
	\gamma_1	& \triangleq  {\frac{8 \log n}{1 - \lambda} \log\frac{1}{\alpha}},   \ \
	\gamma_2\triangleq \frac{1}{n}\, \min_ {\theta\notin\Theta^*}
	\left(F(\theta) -F(\theta^*)\right),
	\end{align*}
	with $\alpha$ from Assumption~\ref{assum:init}(b) and $\lambda = 1 - \frac{1}{18U}$.
\end{theorem} 

The bound of Theorem \ref{thm_lineal} is an improvement by a factor of $n$ compared to 
the bounds of Theorem \ref{teo2}~\eqref{linear_bayes} when the graphs are static. 
Indeed, the term $1/(1-\lambda)$ is $\mathcal{O}(n^2)$ in Theorem \ref{teo2} if $B=1$; the same term is $\mathcal{O}(U)$ in Theorem \ref{thm_lineal} and, assuming $U$ is within a constant factor of $n$, this becomes $\mathcal{O}(n)$. We note, however, that the requirements of this theorem are more stringent than those of Theorem \ref{teo2}. Not only does the graph have to be fixed, but all nodes need to know an upper bound $U$ on the total number of agents.
Moreover, the bound $U$ has to be within a constant factor of the actual number of agents.
More details on fast algorithms for distributed optimization and learning can be found in a tutorial paper~\cite{olsh-cdc2016}.

\subsection{Directed Time-Varying Graphs}

In \cite{ned16}, the authors proposed a new algorithm inspired by the Push-Sum Protocol that is able to guarantee convergence for directed graphs, given as follows:

\begin{subequations}\label{push_bayes}
	\begin{align}
	y^i_{k+1} & =\sum\limits_{j \in N_k^i } \frac{y_k^j} {d^j_k}, \\
	\mu_{k+1}^i\left(\theta\right) & \propto 
	\left(\prod\limits_{j \in N_k^i }\mu_{k}^j\left(\theta\right)^{{y_k^j} / {d^j_k}}\ell^i(s_{k+1}^i|\theta)\right)^{{1} / {y_{k+1}^i}}.
	\end{align}
\end{subequations}

For this algorithm, we have the following result about its convergence behavior.
\begin{theorem}\label{thm_1}
	Assume that the graph sequence $\{\mathcal{G}_k\}$ is $B$-strongly connected and that assumption \ref{assum:init} holds. Also, let $\rho\in(0,1)$ be a given error percentile (or confidence value).
	Then, for the update rule in Eq.~(\ref{push_bayes}), with $y_0^i=1$ and uniform initial beliefs, has the following property:
	there is an integer $\boldsymbol{N}(\rho)$ such that, with probability $1 -\rho$, 
	for all $k\ge \boldsymbol{N}(\rho)$ and for all $\theta\notin\Theta^*$ 
	there holds
	$$\mu_{k}^i\left(\theta\right) \leq \exp\left( -\frac{k}{2} \gamma_2+ \frac{1}{\delta}\gamma_1^i\right)
	\quad\hbox{for all } i = 1, \ldots, n,$$ 
	where $	\boldsymbol{N}(\rho)
	\triangleq \left\lceil\frac{8 \left( \log\left(\alpha \right)\right) ^2  \log\left(\frac{1}{\rho} \right) }{\delta^2\gamma_2^2} + 1 \right\rceil,
	$
	with $\gamma_1^i$ and $\gamma_2$ as defined in Theorem \ref{teo2}.
	The constants $C$, $\delta$ and $\lambda$ satisfy the following relations:\\
	\noindent (1)
	For general $B$-strongly-connected graph sequences $\{\mathcal{G}_k\}$, we have
	\begin{align*}
	C = 4, & \ \ \ \ \lambda = (1-{1} / {n^{nB}})^{\frac{1}{B}}, \ \ \ \delta \geq {1} / {n^{nB}}. 
	\end{align*}
	\noindent (2) If every graph $G_k$ is regular with $B=1$, then we have
	\begin{align*}
	C = \sqrt{2}, & \ \ \ \ \lambda = (1-{1} / {4n^3} )^{\frac{1}{B}}, \ \ \ \delta = 1.
	\end{align*} 
\end{theorem}

\subsection{Infinite Sets of Hypotheses}

All previously discussed results assume that the hypothesis set is finite. 
The exponential convergence rates discussed so far depend on some form of 
distance between the optimal hypothesis and the second best one, and such results are extendable to the case of countably many hypotheses. However, in the case of a continuum of hypothesis, 
this approach will encounter obstacles. 
In a recent work \cite{nedic2016b}, the exponential rates have been established for a compact set $\Theta$ of hypotheses. 
In this case, the update rule for a measurable set $B \subseteq \Theta$ is defined as
\begin{align}\label{protocol_cont}
\mu_{k}^i\left(B\right) & \propto \int\limits_{\theta \in B} \prod\limits_{j=1}^{n} \left( \frac{d\mu_{k-1}^j\left(\theta\right)}{d\lambda(\theta)}\right)^{A_{ij}}  \ell^i(s_{k}^i|\theta) d\lambda\left(\theta\right),
\end{align}
where $\lambda$ is a measure with respect to which every belief $\mu_k^i$ is absolutely continuous.
The particular details of the rate results can be found in~\cite{nedic2016b}.

\section{Conclusions and Future Work}\label{conclusions}

We presented a highlight of recent developments on the problem of distributed (non-Bayesian) learning. We discussed the problem statement and how different assumptions on the learning model, communication graphs and hypothesis sets lead to different algorithmic implementations. 
We showed that the original Bayesian approach can be interpreted as a method for solving a related optimization problem. 

Future work should focus on models where the observations are not necessarily identically distributed nor independent. Recent results on the concentration measures without independence provide the theoretical foundations for getting non-asymptotic rates in more general cases \cite{kon16,duc12}. 
Such time dependence could model changes in the optimal hypotheses, changes on the likelihood models or the Bregman divergences used. Online optimization have been shown efficient for some time dependencies \cite{mok16,sha16}.
\bibliographystyle{IEEEtran} 

\bibliography{IEEEfull,bayes_cons_2}

\end{document}